\documentclass[11pt]{article}
\usepackage{amsfonts,amssymb}
\usepackage{diagrams}

\def\longratto{\longrightarrow} 

\newcounter{paragrafsubsub}[subsubsection]
\renewcommand{\theparagrafsubsub}{%
\thesubsubsection.\roman{paragrafsubsub}}
\newcommand{\paragrafsubsub}{%
\refstepcounter{paragrafsubsub}
{\bf \theparagrafsubsub}\hspace{0.2em}--- }

\newcounter{paragrafsub}[subsection]
\renewcommand{\theparagrafsub}{\thesubsection.\arabic{paragrafsub}}
\newcommand{\paragrafsub}{%
\refstepcounter{paragrafsub}
{\bf \theparagrafsub}\hspace{0.2em}--- }

\newcounter{paragraf}[section]
\renewcommand{\theparagraf}{\thesection.\arabic{paragraf}}
\newcommand{\paragraf}{%
\refstepcounter{paragraf}
{\bf \theparagraf}\hspace{0.2em}--- }

\newcommand\paragraphe{%
\par \indent
\ifcase\value{subsection} %
\paragraf
\else
\ifcase\value{subsubsection}\paragrafsub %
\else\paragrafsubsub
\fi\fi
}

\def\longto{\longrightarrow}

\def\cone{{\mathcal C}}
\def\tc{{\mathcal{TC}}}
\def\ac{{\mathcal{AC}}}
\def\sac{{\mathcal{SAC}}}
\def\lr{{\mathcal{LR}}}
\def\Chi{X}

\def\PP{{\mathbb P}}
\def\QQ{{\mathbb Q}}\def\ZZ{{\mathbb Z}}

\def\kk{{\mathbb K}}

\def\lg{{\mathfrak g}}

\def\hlg{{\hat{\lg}}}

\def\Lim{{p_\lambda}}
\def\kbprod{\odot_0}

\def\Ker{{\rm Ker}}
\def\Hom{{\rm Hom}}
\def\Span{{\rm Span}}
\def\Orb{{\mathcal O}}

\def\LH{L}
\def\Li{{\mathcal{L}}}
\def\Mi{{\mathcal{M}}}

\def\XssL{X^{\rm ss}(\Li)}
\def\CGLX{\ac^G_\Lambda(X)}

\def\quot{/\hspace{-.5ex}/}
\def\Om{\Omega}

\def\Face{{\mathcal F}}

\def\Pic{\rm Pic}

\def\hG{{\hat{G}}}
\def\hB{{\hat{B}}}
\def\hT{{\hat{T}}}
\def\hW{{\hat{W}}}\def\hw{{\hat{w}}}
\def\hnu{{\hat\nu}}
\def\hP{{\hat P}}

\def\P{{\mathcal P}}

\newtheorem{lemma}{Lemma}
\newtheorem{prop}{Proposition}
\newtheorem{theo}{Theorem}

\newtheorem{theoi}{Theorem}

\newenvironment{proof}{{\noindent\bf Proof.}}{\hfill $\square$}
\newenvironment{defin}{{\noindent\bf Definition.}}{\\}
\newenvironment{remark}{{\noindent\bf Remark.}}{}
\newenvironment{nota}{{\noindent\bf Notation.}}{\\}

\begin{document}
\title{Geometric Invariant Theory and\\
Generalized Eigenvalue Problem II}
\author{N. Ressayre}

\maketitle

\begin{abstract}
Let $G$ be a connected reductive subgroup of a complex connected reductive group $\hat{G}$.
Fix maximal tori and Borel subgroups of $G$ and $\hG$.
Consider the cone $\lr^\circ(\hG,G)$ generated by the pairs $(\nu,\hnu)$ of strictly dominant characters 
such that $V_\nu$  is a submodule of $V_{\hnu}$.
The main result of this article is a bijective parametrisation of the faces of $\lr^\circ(\hG,G)$.
We also explain when such a face is contained in another one.

In way, we obtain results about the faces of the Dolgachev-Hu's $G$-ample cone. 
We also apply our results to reprove known results about the moment polytopes.
\end{abstract}

\section{Introduction}
 
Let $G$ be a connected reductive subgroup of a complex connected reductive group $\hat{G}$.
Fix maximal tori and Borel subgroups of $G$ and $\hG$.
Consider the cone $\lr^\circ(\hG,G)$ generated by the pairs $(\nu,\hnu)$ of strictly dominant characters 
such that $V_\nu$  is a submodule of $V_{\hnu}$.
This work is a continuation of \cite{GITEigen}.
We obtain results about general GIT-cones and apply it to obtain a bijective parametrisation of the faces
of the cone  $\lr^\circ(\hG,G)$.

Consider a connected reductive group $G$ acting on a projective variety $X$.
To any $G$-linearized line bundle $\Li$ on $X$ we associate the following open subset 
$X^{\rm ss}(\Li)$  of $X$:
$$
X^{\rm ss}(\Li)=\left \{
x\in X \mbox{ : }
\exists n>0\mbox{ and }
\sigma\in{\rm H}^0(X,\Li^{\otimes n})^G \mbox{ such that }\sigma(x)\neq 0
\right \}.
$$
The points of $X^{\rm ss}(\Li)$ are said to be {\it semistable} for $\Li$.
Note that if $\Li$ is not ample, this notion of semistability is not the standard one. 
In particular, the quotient $\pi_\Li\::\:X^{\rm ss}(\Li)\longto X^{\rm ss}(\Li)\quot G$ is a 
good quotient, if $\Li$ is ample.
In this context, we ask for:
\begin{center}
  What are the $\Li$'s with non empty set $X^{\rm ss}(\Li)$ ?
\end{center}
Let us fix a freely finitely generated subgroup $\Lambda$ of the group $\Pic^G(X)$ of $G$-linearized
line bundles on $X$. 
Let $\Lambda_\QQ$ denote  the $\QQ$-vector space containing $\Lambda$ as a lattice.
Consider the convex cones $\tc_\Lambda^G(X)$ (resp. 
$\ac_\Lambda^G(X)$) generated in $\Lambda_\QQ$ by the $\Li$'s (resp. the ample $\Li$'s) in $\Lambda$ 
which have non zero $G$-invariant sections. By~\cite{DH} (see also \cite{GeomDedic}), 
$\ac^G_\Lambda(X)$ is a closed convex rational polyhedral cone in the dominant cone of $\Lambda_\QQ$.
We are interested in the faces of $\ac^G_\Lambda(X)$ and $\tc_\Lambda^G(X)$.\\

We need to introduce a definition due to D.~Luna.
Assume that $X$ is smooth. Let $\Orb$ be an orbit of $G$ in $X$.
For $x\in\Orb$, we consider the action of the isotropy $G_x$ on the normal space $N_x$ of $\Orb$ 
in $X$ at $x$. The pair $(G_x,N_x)$ is called the {\it type} of the orbit $\Orb$ and is defined up
to conjugacy by $G$. The main part of Theorem~\ref{th:HFCF} is:

\begin{theoi}
\label{thi:typeOrb}
We assume that $X$ is smooth.
Let $\Face$ be a face of $\ac_\Lambda^G(X)$. 

Then, the type  of the closed orbit in  $\pi_\Li^{-1}(\xi)$ for $\xi\in X^{\rm ss}(\Li)\quot G$ 
general does not depends on the choice of an ample $G$-linearized line bundle $\Li$ in the
 relative interior of $\Face$.

We will call this type the type of $\Face$.
\end{theoi}

Let $\Face$ be a face of $\ac^G(X)$. Let $\Li_0$ be any point in the relative interior of $\Face$. 
The local geometry $\ac^G(X)$ around  $\Face$ is described by the convex cone $\cone_\Face$ 
generated by the vectors $p-\Li_0$ for $p\in\ac^G(X)$. 

We now introduce some notation to describe this cone. 
Consider the quotient $\pi\,:\, X^{\rm ss}(\Li)\longto X^{\rm ss}(\Li)\quot G$. 
Let $x$ be any point in $X^{\rm ss}(\Li)$ with closed orbit  in  $X^{\rm ss}(\Li)$, and so, reductive 
isotropy $G_x$.
Then, the fiber $\pi^{-1}(\pi(x))$ is isomorphic to a fiber product $G\times_{G_x}L$, for an affine
$G_x$-variety $L$ with a fixed point as unique closed orbit. 
Let $X(G_x)$ denote the group of characters of $G_x$. 
Consider the semicone $\cone_x$ in $X(G_x)\otimes \QQ$ generated by the weights of $G_x$ on the
set of regular functions on $L$.
Finally, we consider is  the linear map $\mu\,:\,\Lambda_\QQ\longto X(G_x)\otimes \QQ$ obtained by considering 
the action of $G_x$ on the fibers $\Li_x$ in $\Li\in\Lambda$ over $x$.
 
\begin{theoi}
\label{thi:coneloc}
With above notation, we have:
$$
\cone_\Face=\mu^{-1}(\cone_x).
$$ 
\end{theoi}

In the symplectic setting, S.~Sjamaar obtained a description of the local
structure of the moment polytope (see~\cite{Sja:convreex}) which is closed from
Theorem~\ref{thi:coneloc}. In Sjamaar's situation, $G_x^\circ$ is a torus 
which simplifies a little bit.

Now, assume that the variety $X$ equals $Y\times G/B$, for a $G$-variety $Y$.
Let $\Li$ be an ample $G$-linearized line bundle on $Y$.
Let $\Lambda$ be the subgroup of $\Pic^G(X)$ generated by the pullback of $\Li$ and the
pullbacks of the $G$-linearized line bundles on $G/B$.
Then, $\tc^G_\Lambda(X)$ is a cone over the moment polytope $P(Y,\Li)$ 
defined in \cite{Br:genface};
in particular, the faces of $\tc^G_\Lambda(X)$ correspond bijectively to the faces of $P(Y,\Li)$.

Following \cite{Br:genface}, we show in Proposition~\ref{prop:PXLempty} below, that any moment polytope
$P(Y,\Li)$ can be describe in terms of one which intersects the interior of the dominant chamber. 
We now assume that $P(Y,\Li)$ intersects the interior of the dominant chamber and that $Y$ is smooth.
In Proposition~\ref{prop:PYL} below, we associate to each face of $P(Y,\Li)$ which intersects the interior
of the dominant chamber a well $B$-covering pair (see Definition~\ref{def:wellBcov}) of $Y$ improving 
(with stronger assumptions) \cite[Theorem~1 and 2]{Br:genface}.\\

Now, $G$ is assumed to be embedded in another connected reductive group $\hG$.
We fix maximal tori $T\subset\hT$ and Borel subgroups $B\supset T$ and $\hB\supset \hT$ of $G$ and $\hG$.
Consider the diagonal action of $G$ on $\hG/\hB\times G/B$ and the associated GIT-cone
$\ac^G(\hG/\hB\times G/B)$. 
Actually, $\ac^G(\hG/\hB\times G/B)$ identifies with the cone generated by pairs 
$(\nu,\hnu)$ of strictly dominant character of $T\times\hT$ such that the dual of the $G$-module 
associated to $\nu$ can be $G$-equivariantly 
embedded in the $\hG$-module associated to $\hnu$.
The interior of $\ac^G(\hG/\hB\times G/B)$ is non empty if and only if no non trivial 
connected normal subgroup
of $G$ is normal in $\hG$: we assume, from now on that $\ac^G(\hG/\hB\times G/B)$ has non empty interior.
Theorem~\ref{th:ppal} below gives a bijective parametrisation of the faces of $\ac^G(\hG/\hB\times G/B)$.
Moreover, we can read very easily the inclusions between faces using this parametrisation. 
To avoid too many notation, in this introduction, we will only state our results in the case when $\hG=G^2$.

For any standard parabolic subgroup $P$ of $G$, we consider the cohomology group ${\rm H}^*(G/P,\ZZ)$
and its basis consisting in classes of Schubert varieties. 
We consider on this group the Belkale-Kumar product $\kbprod$ defined in \cite{BK}.
The coefficient-structure of this product in this basis are either $0$ or the  coefficient-structure of
the usual cup product. These coefficients are parametrized by the triple of Schubert classes.

\begin{theoi}
The group $G$ is assumed to be semi-simple.
  The set of faces of  $\ac^G((G/B)^3)$ correspond bijectively to the set of structure coefficient of 
$({\rm H}^*(G/P,\ZZ),\kbprod)$ equal to one, for the various standard parabolic subgroups $P$ of $G$. 
\end{theoi}

We will now explain how to read the inclusion off this parametrization.
Let $P$ and $P'$ be two standard parabolic subgroups. 
Let $\Lambda_1,\Lambda_2$ and $\Lambda_3$ (resp. $\Lambda_1',\Lambda_2'$ and $\Lambda_3'$) three Schubert varieties in 
$G/P$ (resp. $G/P'$) such the corresponding coefficients structure for $\kbprod$ 
equal to one. Let $\Face$ and $\Face'$ denote the corresponding faces of $\ac^G((G/B)^3)$.

\begin{theoi}
Let $\Face$ and $\Face'$ be two faces of  $\ac^G((G/B)^3)$.
The following are equivalent:
\begin{enumerate}
\item $\Face\subset\Face'$;
\item $P\subset P'$ and $\pi(\Lambda_i)=\Lambda'_i$ for $i=1,2$ and $3$, where $\pi\,:\,G/P\longto G/P'$ 
is the $G$-equivariant morphism mapping $P/P$ on $P'/P'$.
\end{enumerate}
\end{theoi}

\noindent{\bf Convention.}
The ground field $\kk$ is assumed to be algebraically closed of characteristic zero.
The notation introduced in the environments ``{\bf Notation.}'' are fixed for all the sequence 
of the article.

\section{An example of GIT-cone}
\label{sec:Cov}

Let us fix a connected reductive group $G$ acting on an irreducible projective 
algebraic variety $X$.

\subsection{An Ad Hoc notion of semistability}

As in the introduction, for any  $G$-linearized line bundle $\Li$ on $X$, we consider the
following set of {\it semistable points}:
$$
X^{\rm ss}(\Li)=\left \{
x\in X \mbox{ : }
\exists n>0\mbox{ and }
\sigma\in{\rm H}^0(X,\Li^{\otimes n})^G \mbox{ such that }\sigma(x)\neq 0
\right \}.
$$ 
To precise the acting group, we sometimes denote $X^{\rm ss}(\Li)$ by $X^{\rm ss}(\Li,G)$.

The subset $X^{\rm ss}(\Li)$ is open and stable by $G$.
  Note that this definition of $X^{\rm ss}(\Li)$ is  the standard one only when $\Li$ is
ample. Indeed, one usually imposes that the open
subset defined by the non vanishing of $\sigma$ to be affine.

If $\Li$ is ample, there exists a categorical quotient:
$$
\pi\::\:X^{\rm ss}(\Li)\longto X^{\rm ss}(\Li)\quot G,
$$
such that $X^{\rm ss}(\Li)\quot G$ is a projective variety and $\pi$ is affine.
A point $x\in X^{\rm ss}(\Li)$ is said to be {\it stable} if $G_x$ is finite and $G.x$ is 
closed in $X^{\rm ss}(\Li)$. 
Then, for all stable point $x$ we have $\pi^{-1}(\pi(x))=G.x$; and the set 
$X^{\rm s}(\Li)$ of stable points is open in $X$.

\subsection{Definitions}

Let us recall from the introduction that $\Lambda$ is a freely finitely generated subgroup of $\Pic^G(X)$ 
and $\Lambda_\QQ$ is the $\QQ$-vector space containing $\Lambda$ as a lattice.
Since $X^{\rm ss}(\Li)=X^{\rm ss}(\Li^{\otimes n})$, for any $G$-linearized line bundle and
any positive integer $n$, we can define $X^{\rm ss}(\Li)$ for any $\Li\in\Lambda_\QQ$.
We consider the following {\it total $G$-cone}:
$$
\tc^G_\Lambda(X)=\{\Li\in\Lambda_\QQ\;:\;  X^{\rm ss}(\Li) {\rm\ is\ not\ empty}\}.
$$
Since the tensor product of two non zero $G$-invariant sections is a non zero $G$-invariant section, 
$\tc^G_\Lambda(X)$ is a convex cone.

Consider the convex cones   $\Lambda^{+}_\QQ$ and $\Lambda^{++}_\QQ$ generated respectively by the semiample 
and ample elements of $\Lambda$. 
For all $\Li\in\Lambda^+_\QQ$ (resp. $\Lambda^{++}_\QQ$), there exists a positive integer $n$ such that
$\Li^{\otimes n}$ is a semiample (resp. ample) $G$-linearized line bundle on $X$ in $\Lambda$.
So, any set of semistable points associated to a point in $\Lambda^{+}_\QQ$ (resp. $\Lambda^{++}_\QQ$) 
is in fact a set  of semistable point associated to  a semiample (resp. ample) $G$-linearized line bundle.
We consider the following {\it semiample and ample $G$-cones}:
$$
\sac^G_\Lambda(X)=\tc^G_\Lambda(X)\cap\Lambda^{+}_\QQ {\rm\ \ and\ \ }   \ac^G_\Lambda(X)=\tc^G_\Lambda(X)\cap\Lambda^{++}_\QQ.
$$
By~\cite{DH} (see also \cite{GeomDedic}), $\ac^G_\Lambda(X)$ is a closed convex rational polyhedral cone in $\Lambda^{++}_\QQ$. This cone is the central object of this article.

Two points $\Li$ and $\Li'$ in $\CGLX$ are said to be {\it GIT-equivalent} if $X^{\rm ss}(\Li)=X^{\rm ss}(\Li')$.
An equivalence class is simply called a {\it GIT-class}. 
 
For $x\in X$, the {\it stability set of $x$} is the set of $\Li\in\Lambda_\QQ^{++}$ such that 
$X^{\rm ss}(\Li)$ contains $x$; it is denoted by $\Omega_\Lambda(x)$ 
or  $\Omega_\Lambda(G.x)$.
In \cite{GeomDedic}, we have studied the geometry of the GIT-classes and the stability sets with lightly different 
assumptions (no $\Lambda$ for example). However all the results and proofs of \cite{GeomDedic} remain valuable here.
In particular, there are only finitely many GIT-classes; and
each GIT-class is the relative interior of a closed convex polyhedral cone of $\Lambda_\QQ^{++}$.

\subsection{An example of $G$-ample cone}
\label{sec:exple}

\begin{nota}
If $\Gamma$ is an affine algebraic group, $[\Gamma,\Gamma]$ will denote its derived subgroup
and $X(\Gamma)$ will denote the character group of $\Gamma$.  
\end{nota}

For later use, we consider here a $G$-ample cone for the action of $G$ over an affine variety.
More precisely, let $V$ be an affine $G$-variety containing a fix point $O$ as unique closed orbit. 
The action of $G$ over the fiber gives a morphism $\mu^\bullet(O,G)\ :\ \Pic^G(V)\longto \Chi(G)$.
By \cite[Lemma~7]{GeomDedic}, $\mu^\bullet(O,G)$ is an isomorphism.
We denote by $V^{\rm ss}(\chi)$ the set of semistable points for the $G$-linearized line bundle $\Li_\chi$ 
associated to $\chi\in\Chi(G)$; that is, the trivial line bundle on $V$ linearized by $\chi$.
As in the projective case, we consider the $G$-ample cone $\ac^G(V)$ in $\Chi(G)\otimes\QQ$.

For any $\chi\in \Chi(G)$, we have:
$$
{\rm H}^0(V,\Li_\chi)^G=\{f\in\kk[V]\ :\ \forall x\in V\ \ (g.f)(x)=\chi(g)f(x)\}=\kk[V]_\chi.
$$
Note that ${\rm H}^0(V,\Li_\chi)^G$ is contained in $\kk[V]^{[G,G]}$. Set
$$
S=\{\chi\in\Chi(G)\ :\ {\rm H}^0(V,\Li_\chi)^G\ {\rm\ is\ non\ trivial}\}.
$$
It is the set of weights of $G/[G,G]$ in $\kk[V]^{[G,G]}$.
We have:

\begin{lemma}
\label{lem:CGV}
We assume that $V$ is irreducible.
The set $S$ is a finitely generated semigroup in $\Chi(G)$.
Moreover, $\ac^G(V)$ is the convex cone generated by $S$; 
it is strictly convex.
\end{lemma}

\begin{proof}
Since $\kk[V]^{[G,G]}$ is a finitely generated  algebra, $S$ is a finitely generated semigroup. 
The fact that $\ac^G(V)$ is generated by $S$ is obvious.
Finally,  $\ac^G(V)$ is strictly convex since ${\rm H}^0(V,\Li_0)^G=\kk$.
\end{proof}

\section{Slice Etale Theorem}
\label{sec:Luna}

In this section, we recall some very useful results of D.~Luna.
We fix an ample $G$-linearized line bundle $\Li$ on the irreducible projective $G$-variety $X$.

\subsection{Closed orbits in general position}

\begin{nota}
If $H$ is a subgroup of $G$, $N_G(H)$ denotes the normalizer of $H$ in $G$.
Consider the quotient $\pi\::\: X^{\rm ss}(\Li)\longto X^{\rm ss}(\Li)\quot G$.
For all $\xi\in X^{\rm ss}(\Li)\quot G$, 
we denote by $T(\xi)$ the unique closed orbit of $G$ in $\pi^{-1}(\xi)$.
We denote by $(X^{\rm ss}(\Li)\quot G)_{\rm pr}$ the set of $\xi$ such that there
exists an open neighborhood $\Omega$ of $\xi$ in $X^{\rm ss}(\Li)\quot G$ such that 
the orbits $T(\xi')$ are isomorphic to $T(\xi)$, for all
$\xi'\in \Omega$.
\end{nota}

Since $\pi$ is a gluing of affine quotients, some results on the actions 
of $G$ on affine variety remains true for $X^{\rm ss}(\Li)$. For example, the following theorem 
is a result of Luna and Richardson (see \cite[Section~3]{LunaRichardson:chevalley} and
\cite[Corollary~4]{Luna:adh} or \cite[Section~7]{PV}):

\begin{theo}
\label{th:cogp}
With above notation, if $X$ is normal, we have:
\begin{enumerate}
\item \label{ass:genisot}
The set $(X^{\rm ss}(\Li)\quot G)_{\rm pr}$ is a non empty open subset of $X^{\rm ss}(\Li)\quot G$.
Let $H$ be the isotropy of a point in $T(\xi_0)$ with $\xi_0\in (X^{\rm ss}(\Li)\quot G)_{\rm pr}$.
The group $H$ has fixed points in $T(\xi)$ for any $\xi\in X^{\rm ss}(\Li)\quot G$.

\item Let $Y$ be the closure of $\pi^{-1} \left((X^{\rm ss}(\Li)\quot G)_{\rm pr}\right)^H$ in $X$.
It is the union of some components of $X^H$.
Then, $H$ acts trivially on some positive power $\Li_{|Y}^{\otimes n}$ of $\Li_{|Y}$. Moreover, the natural map
$$
Y^{\rm ss}(\Li_{|Y}^{\otimes n})\quot (N_G(H)/H)\longto X^{\rm ss}(\Li)\quot G
$$
is an isomorphism. Moreover, $Y$ contains stable points for the action of $N_G(H)/H$ and the line bundle
$\Li_{|Y}^{\otimes n}$.
\end{enumerate}
\end{theo}

A subgroup $H$ as in Theorem~\ref{th:cogp} will be called {\it 
a generic closed isotropy of $X^{\rm ss}(\Li)$}.
The conjugacy class of $H$ which is obviously unique is called 
{\it the generic closed isotropy of $X^{\rm ss}(\Li)$}.

\subsection{The principal Luna stratum}
\label{sec:ppalstrate}

When $X$ is smooth, the open subset $(\XssL\quot G)_{\rm pr}$ is called the
{\it principal Lana stratum} and has very useful properties 
(see~\cite{Luna:slice} or \cite{PV}): 

\begin{theo}[Luna]
\label{th:slice}
We assume that $X$ is smooth.
Let $H$ be a generic closed isotropy of $X^{\rm ss}(\Li)$.

Then, there exists  a $H$-module $\LH$ such that
for any $\xi\in(\XssL\quot G)_{\rm pr}$, there exist points $x$ in $T(\xi)$  satisfying: 
\begin{enumerate}
\item $G_x=H$;
\item the $H$-module $T_xX/T_x(G.x)$ is isomorphic to the sum of $L$ and its fix points,
in particular, it is independent of $\xi$ and $x$; 
\item for any $v\in L$, $0$ belongs to the closure of $H.v$;
\item the fiber $\pi^{-1}(\xi)$ is isomorphic to $G\times_HL$. 
\end{enumerate}
\end{theo}

\subsection{The fibers of quotient morphisms}

Another useful consequence of Luna's Slice Etale Theorem is (see~\cite{Luna:slice} or \cite{PV}):

\begin{prop}
\label{prop:fibrepi}
  Let $x$ be a semistable point for $\Li$ whose the orbit is closed in $X^{\rm ss}(\Li)$.
Then, there exists an affine $G_x$-variety $V$ containing a unique closed orbit which is a fixed point and such 
that $\pi^{-1}(\pi(x))$ is isomorphic to $G\times_{G_x}V$.
\end{prop}

\section{About faces of the $G$-ample cone}
\label{sec:gen}

\subsection{Isotropy subgroups associated to faces of $\CGLX$}

Let $\varphi$ be a linear form on $\Lambda_\QQ$ which is non negative on $\CGLX$.
If the set of $\Li\in\CGLX$ such that $\varphi(\Li)=0$ is non empty it will be called 
{\it a face of $\CGLX$}.
Now, we associate two invariants to a face $\Face$ of $\CGLX$.

\begin{theo}
\label{th:HFCF}
Let $\Face$ be a face of $\CGLX$. 
Let $\Li$ be  a point in the relative interior of $\Face$. 
Then, we have:

\begin{enumerate}
\item \label{ass:HF}
The generic closed isotropy of $X^{\rm ss}(\Li)$ does not depend on the point
$\Li$ in the relative interior of $\Face$, but only in $\Face$.
We call this isotropy the {\it generic closed isotropy of $\Face$}.

Let us fix a generic closed isotropy $H$ of $\Face$. 
\item \label{ass:gciL}
For any $\Mi\in\Face$, $H$ fixes points in any closed orbit of $G$ in $X^{\rm ss}(\Mi)$.

\item \label{ass:YF} The closure $Y$ of 
$\bigg( \pi_\Li^{-1} \left((X^{\rm ss}(\Li)\quot G)_{\rm pr}\right)\bigg)^H$ in $X$
does not depends on a choice of $\Li$.
Let $Y_\Face$ denote this subvariety of $X^H$; it is the union of some components of $X^H$.

\item \label{ass:YL}
The group  $H$ acts trivially on some positive power $\Li_{|Y_\Face}^{\otimes n}$ of $\Li_{|Y_\Face}$. Moreover, the natural map
$$
Y_\Face^{\rm ss}(\Li_{|Y_\Face}^{\otimes n})\quot (N_G(H)/H)\longto X^{\rm ss}(\Li)\quot G
$$
is an isomorphism.
Moreover, 
$Y_\Face$ contains stable points for the action of $N_G(H)/H$ and the line bundle $\Li_{|Y_\Face}^{\otimes n}$.

\item \label{ass:Y+}
Set $Y^+_\Face:=\{x\in X\::\:\overline{H.x}\cap Y_\Face\neq\emptyset\}$.
Then $G.Y^+_\Face$ contains an open subset of $X$.

\end{enumerate}
\end{theo}

\begin{proof}
Let $\Li_1,\,\Li_2\in \CGLX$.
By an easy argument of convexity, to prove Assertion~\ref{ass:HF} 
it is sufficient to prove that the generic closed isotropy of $X^{\rm ss}(\Li)$
does not depend on $\Li$ in the open segment $]\Li_1,\Li_2[$.
Let us fix $\Li,\,\Li'\in ]\Li_1,\Li_2[$.
Let $x\in X$ which maps in $(X^{\rm ss}(\Mi)\quot G)_{\rm pr}$, for 
$\Mi=\Li_1,\,\Li_2,\,\Li$ and $\Li'$ by the quotient maps.

Recall that $\Om_\Lambda(x)$ is a polyhedral convex cone.
Since  $\Li_1$ and $\Li_2$ belong to $\Om_\Lambda(x)$,
$\Li$ and $\Li'$ belong to the relative interior of the same face of 
$\Om_\Lambda(x)$.
By  \cite[Proposition~6, Assertion~(iii)]{GeomDedic}, there exists $x'\in\overline{G.x}$
such that this face is  $\Om_\Lambda(x')$. 
But, \cite[Proposition~6, Assertion~(i)]{GeomDedic} shows that 
the closed orbits in $X^{\rm ss}(\Li)\cap \overline{G.x'}$ and 
$X^{\rm ss}(\Li')\cap \overline{G.x'}$ are equal.
Now, our choice of the point $x$ implies that the generic closed isotropies of 
$X^{\rm ss}(\Li)$ and $X^{\rm ss}(\Li')$ are equal.\\

Let $H$ be a generic closed isotropy of  $X^{\rm ss}(\Li)$.
Let $Y$ be the closure of 
$X^H\cap\pi_\Li^{-1}\left((X^{\rm ss}(\Li)\quot G)_{\rm pr}\right)$.
By Theorem~\ref{th:cogp}, $N_G(H)$ acts transitively on the set of irreducible components of $Y$.
Let $Y_1$ be such a component of $X^H$.
Again by Theorem~\ref{th:cogp}, $\pi_\Li(Y_1\cap X^{\rm ss}(\Li))=X^{\rm ss}(\Li)\quot G$;
that is, any closed $G$-orbit in $X^{\rm ss}(\Li)$ intersects $Y_1$.
Finally, $Y$ is the union of irreducible components of $X^H$ which intersect a general 
closed $G$-orbit in $X^{\rm ss}(\Li)$. 
But, the above proof of Assertion~\ref{ass:HF} shows that a general closed orbit in 
$X^{\rm ss}(\Li)$ is also a closed orbit in $X^{\rm ss}(\Li')$.
In particular, $Y$ is the closure of $X^H\cap\pi_{\Li'}^{-1}\left((X^{\rm ss}(\Li')\quot G)_{\rm pr}\right)$.
Assertion~\ref{ass:YF} follows.\\

Let us now fix a generic closed isotropy $H$ of $\Face$.
Let $\Mi_1\in \Face$.
By Assertion~\ref{ass:genisot} of Theorem~\ref{th:cogp}, to prove the 
second assertion, it is sufficient to prove that the generic closed isotropy
of $X^{\rm ss}(\Mi_1)$ contains $H$.
By \cite[Theorem~4]{GeomDedic}, there exists a point $\Mi_2$ in the relative 
interior of $\Face$ such that $X^{\rm ss}(\Mi_1)$ contains $X^{\rm ss}(\Mi_2)$.
The inclusion $X^{\rm ss}(\Mi_2)\subset X^{\rm ss}(\Mi_1)$ induces a surjective
morphism $\eta\,:\,X^{\rm ss}(\Mi_2)\quot G\longto X^{\rm ss}(\Mi_1)\quot G$.
Let $\xi'\in (X^{\rm ss}(\Mi_2)\quot G)_{\rm pr}$ such that 
$\xi=\eta(\xi')\in (X^{\rm ss}(\Mi_1)\quot G)_{\rm pr}$.
Let $x$ be a point in the closed $G$-orbit in $X^{\rm ss}(\Mi_1)$ over $\xi$.
The fiber in $X^{\rm ss}(\Mi_1)$ over $\xi$ is fibered over $G.x$; hence, 
for any $y$ in this fiber, $G_y$ is conjugated to a subgroup of $G_x$.
Since this fiber contains the fiber in $X^{\rm ss}(\Mi_2)$ over 
$\xi'$, $H$ is conjugated to a subgroup of $G_x$.
The second assertion is proved.\\

From now on, $\Li$ is a point in the relative interior of $\Face$.
Let $Y$ be the subvariety of $X^H$ of Assertion~\ref{ass:YF}. 
By Theorem~\ref{th:cogp}, $Y$ satisfies Assertion~\ref{ass:YL}.
Moreover, $G.Y^+$ contains $\pi_\Li^{-1} \left((X^{\rm ss}(\Li)\quot G)_{\rm pr}\right)$; 
and,  Assertion~\ref{ass:Y+} is proved.
\end{proof}

\subsection{Local structure of the $G$-ample cone around a face}

\begin{nota}
  Let $E$ be a prime Cartier divisor on a variety $X$ endowed with a line bundle $\Li$ 
and $\sigma$ be a rational section for $\Li$.
We will denote by $\nu_E(\sigma)\in\ZZ$ the  order of zero of $\sigma$ along $E$.

Let $\P$ be a polyhedron in a rational vector space $V$ and $\Face$ be a face of $\P$.
Let $x$ in the relative interior of $\Face$.
The cone of $V$ generated by the vectors $y-x$ for $y\in\P$ does not depends on the choice of $x$ 
in the relative interior of $\Face$ and will be called {\it the cone of $\P$ viewed from $\Face$}.
This cone carries the geometry of $\P$ in a neighborhood of $x$. 
\end{nota}

Let $\Face$ be a face of $\CGLX$. Let $\Li$ in the relative interior of $\Face$.
Let $x$ be a  semistable point for $\Li$ whose the orbit is closed in $X^{\rm ss}(\Li)$.
Let $V$ be an affine $G_x$-variety satisfying Proposition~\ref{prop:fibrepi}.
Consider the cone $\ac^{G_x}(V)$ as in Section~\ref{sec:exple}.
Notice that $V$ is not necessarily irreducible, and so  $\cone^{G_x}(V)$ is not necessarily convex.

The action of $G_x$ on the fiber over $x$ defines a morphism
$\mu^\bullet(x,G_x)\,:\,\Lambda\longto \Chi(G_x)$; and so, a linear map 
from $\Lambda_\QQ$ on $\Chi(G_x)_\QQ$ also denoted by $\mu^\bullet(x,G_x)$.

\begin{theo}
\label{th:CGXloc}
With above notation, the cone of $\ac^G(X)$ viewed from $\Face$ equals 
$(\mu^\bullet(x,G_x))^{-1}(\ac^{G_x}(V))$.
In particular, if $\mu^\bullet(x,G_x)$ is surjective then $\ac^{G_x}(V)$ is convex.
\end{theo}

\begin{proof}
  Let $\Li_0$ and $\Li$ be two different ample $G$-linearized line bundles in $\Lambda$. 
We assume that $\Li_0$ is the only point in the segment $[\Li;\Li_0]$ which belongs to $\CGLX$.
For convenience, we set $U=X^{\rm ss}(\Li_0)$.
By assumption, there is no $G$-invariant rational section of $\Li$ which is regular on $X$;
we claim that there is no such section which is  regular on $U$.
 
Let us prove the claim. Let us fix a non zero regular $G$-invariant  section $\sigma_0$ of $\Li_0^{\otimes m}$ 
for a positive integer $m$.
Let $\sigma$ be a $G$-invariant rational section of $\Li$ which is regular on $U$.
For any positive integer $k$, $\sigma\otimes\sigma_0^{\otimes k}$ is a rational $G$-invariant section
of $\Li\otimes\Li_0^{\otimes mk}$ which is regular on $U$.
Let $E$ be an irreducible component of codimension one of $X-U$. 
By definition of $U$, $\sigma_0$ is zero along $E$; and, $\nu_E(\sigma_0)>0$.
Then, 
$\nu_E(\sigma\otimes\sigma_0^{\otimes k})=\nu_E(\sigma)+k.\nu_E(\sigma_0)$ is positive for $k$ big enough.
We deduce that $\sigma\otimes\sigma_0^{\otimes k}$ is regular on $X$ for $k$ big enough.
Since by assumption $\Li\otimes\Li_0^{\otimes mk}$ does not belong to $\CGLX$, this implies that 
$\sigma\otimes\sigma_0^{\otimes k}$ and finally $\sigma$ are zero.
The claim is proved.\\

We now fix a point $\Li_0$ in the relative interior of $\Face$ as in the statement.
By an elementary argument of convexity, there exists  an open neighborhood $\Omega$ of $\Li_0$ in 
$\Lambda_\QQ$ such that 
\begin{enumerate}
\item for any $\Li\in \Omega$, if $\Li$ does not belong to $\CGLX$ then
$\Li_0$ is the only point in $[\Li,\Li_0]\cap\CGLX$.
\end{enumerate}
By \cite[Proposition~2.3]{GeomDedic}, we may also assume that for all $\Li\in\Omega$, 
$X^{\rm ss}(\Li)$ is contained in $X^{\rm ss}(\Li_0)$.
It remains to prove that for any $\Li\in \Omega$, $\Li\in\CGLX$ if and only if 
$\mu^\Li(x,G_x)\in \ac^{G_x}(V)$.\\

Let $\Li\in \Omega$ which does not belongs to $\CGLX$. Set $\xi=\pi_{\Li_0}(x)$. 
By the beginning of the proof, for any positive $n$, 
${\rm H}^0(U,\Li^{\otimes n})^G=\{0\}$.
Since $\pi_{\Li_0}^{-1}(\xi)$ is closed in $U$, this implies that 
${\rm H}^0(\pi_{\Li_0}^{-1}(\xi),\Li^{\otimes n})^G=\{0\}$ for all positive $n$.
So, for all positive $n$, ${\rm H}^0(V,\Li_{|V}^{\otimes n})^{G_x}=\{0\}$; and, so
$\mu^\Li(x,G_x)$ does not belong to $\ac^{G_x}(V)$.\\

Let now $\Li\in\Omega\cap\CGLX$.
Since the map $\phi\,:\,X^{\rm ss}(\Li)\quot G\longto X^{\rm ss}(\Li_0)\quot G$ induced by the inclusion
$X^{\rm ss}(\Li)\subset X^{\rm ss}(\Li_0)$ is surjective, there exists $y\in X^{\rm ss}(\Li)$ such that 
$\phi\circ\pi_\Li(y)=\xi$.
Up to changing $y$ by $g.y$ for some $g\in G$, one may assume that $y\in V$.
Let $\sigma$ be a $G$-invariant section of $\Li$ which is non zero at $y$. 
Obviously, the restriction of $\sigma$ is a $G_x$-invariant section of $\Li_{|V}$ 
which is non zero.
It follows that  $\mu^\Li(x,G_x)$ belongs to
$\ac^{G_x}(V)$.

The last assertion follows from  an obvious argue of convexity.
\end{proof}\\

\section{Well covering pairs}

\subsection{The functions $\mu^\bullet(x,\lambda)$}
\label{sec:mu}

Let $\Li\in \Pic^G(X)$.
Let $x$ be a point in $X$ and $\lambda$ be a one parameter subgroup of
$G$.
Since $X$ is complete,  $\lim_{t \to 0}\lambda(t)x$ exists; let $z$
denote this limit.
The image of $\lambda$ fixes $z$ and so the group $\kk^*$ acts
via $\lambda$ on the fiber $\Li_{z}$.
This action defines a character of $\kk^*$, that is, an element of
$\ZZ$ denoted by $\mu^\Li(x,\lambda)$.

The numbers $\mu^\Li(x,\lambda)$ are used in \cite{GIT} to give a 
numerical criterion for stability with respect to an ample $G-$linearized
line bundle $\Li$:
$$
\begin{array}{l}
  x\in X^{\rm ss}(\Li)\iff \mu^\Li(x,\lambda)\leq 0 \mbox{ for all
    one parameter subgroup $\lambda$},\\
x\in X^{\rm s}(\Li)\iff \mu^\Li(x,\lambda)< 0 \mbox{ for all non trivial 
$\lambda$}.
\end{array}
$$

\subsection{Definition}
\label{sec:BB}

\begin{nota}
The set of fix points of the image of $\lambda$ will be denoted by $X^\lambda$; the centralizer of this image
will be denoted by $G^\lambda$. We consider the following parabolic subgroup of $G$: 
$$
P(\lambda)=\left \{
g\in G \::\:
\lim_{t\to 0}\lambda(t).g.\lambda(t)^{-1} 
\mbox{  exists in } G \right \}.
$$
\end{nota}

Let $C$ be an irreducible component of $X^\lambda$. 
Since $G^\lambda$ is connected, $C$ is a $G^\lambda$-stable closed subvariety of $X$.
We set:
$$
C^+:=\{x\in X\::\: \lim_{t\to 0}\lambda(t)x\in C\}.
$$
Then, $C^+$ is a locally closed subvariety of $X$ stable by $P(\lambda)$.
Moreover, the map $\Lim\::\:C^+\longto C,\,x\longmapsto\lim_{t\to 0}\lambda(t)x$ is 
a morphism satisfying:
$$
\forall (l,u)\in G^\lambda\times U(\lambda)\hspace{1cm}
\Lim(lu.x)=l\Lim(x).
$$

Consider over $G\times C^+$ the action of $G\times P(\lambda)$ given by the formula (with obvious notation):
$$
(g,p).(g',y)=(gg'p^{-1},py).
$$
Since the quotient map $G\longto G/P(\lambda)$ is a Zariski-locally trivial principal $P(\lambda)$-bundle; 
one can easily construct a quotient $G\times_{P(\lambda)}C^+$ of $G\times C^+$ by the action 
of $\{e\}\times P(\lambda)$.
The action of $G\times\{e\}$ induces an action of $G$ on $G\times_{P(\lambda)}C^+$.

\begin{defin}
Consider the following $G$-equivariant map
$$
\begin{array}{cccc}
\eta\ :&G\times_{P(\lambda)}C^+&\longto& X\\
       &[g:x]&\longmapsto&g.x.
\end{array}
$$ 
The pair $(C,\lambda)$ is said to be {\it covering} (resp. {\it dominant}) if $\eta$ is birational
(resp. dominant).
It is said to be {\it well covering} if $\eta$ induces an isomorphism from 
$G\times_{P(\lambda)}\Om$ onto a $P(\lambda)$-stable open subset of $X$ 
for an open subset $\Om$ of $C^+$ intersecting $C$.
\end{defin}

\subsection{Face associated to $(C,\lambda)$}

Let us denote by $\mu^\bullet(C,\lambda)$,  the common value of the 
$\mu^\bullet(x,\lambda)$, for $x\in C^+$.

\begin{lemma}
\label{lem:FaceinterC}
  Let $(C,\lambda)$ be a dominant pair.
The set of $\Li\in\CGLX$ such that $\mu^\Li(C,\lambda)=0$ is either empty or a face $\Face$ of $\CGLX$.
Moreover, $\Face$ is the set of $\Li\in\CGLX$ such that $X^{\rm ss}(\Li)$ intersects $C$.

From now on, $\Face$ which only depends on $C$ is denoted by $\Face(C)$.
\end{lemma}

\begin{proof}
  The first assertion is \cite[Lemma~7]{GITEigen}.
If $\Li\in \Face$, then there exists $x\in C^+$ semistable for $\Li$.
By \cite[Lemma~4]{GITEigen}, $\lim_{t\to 0}\lambda(t)x$ is semistable and belongs to $C$.
Conversely, assume that $X^{\rm ss}(\Li)\cap C$ contains $z$.
Since $z$ is fixed by $\lambda$, $\mu^\Li(z,-\lambda)=-\mu^\Li(z,\lambda)$. 
But $z$ is semistable, so $\mu^\Li(z,-\lambda)=0$.
\end{proof}

\section{The case $X=Y\times G/B$}
\label{sec:PYL}

In this section, we assume that $X=Y\times G/B$, 
with a normal projective $G$-variety $Y$.
Moreover, we assume that $\Lambda$ is abundant (see \cite{GITEigen}).

\subsection{General closed isotropy and well covering pairs}

Let $S$ be a torus contained in $G$.
Let $C$ be an irreducible component of $X^S$.

\begin{defin}
The pair $(C,S)$ is said to be {\it admissible} if there exists $x\in C$ such that $G_x^\circ=S$.
The pair is said to be {\it well covering} if there exists a one parameter subgroup $\lambda$  of $S$,
such that $C$ is an irreducible component of $X^\lambda$ and $(C,\lambda)$ is well covering.
\end{defin}

A rephrasing of \cite[Corollary~3]{GITEigen} is

\begin{prop}
\label{prop:ttFC}
We assume that   $X=Y\times G/B$ with a normal projective $G$-variety $Y$ and $\Lambda$ is abundant
(see~\cite[Section~7.4]{GITEigen}).
Let $\Face$ be a face of $\CGLX$ of codimension $r$. Then, there exists an admissible well covering pair
$(C,S)$ with  $S$ of dimension $r$ such that $\Face=\Face(C)$.
\end{prop}

We are now interested in the generic closed isotropy of faces of $\CGLX$:

\begin{prop}
\label{prop:gciGB}
Let $\Face$ be a face of codimension codimension $r$.
Let $(C,S)$ be  an admissible well covering pair with  a $r$-dimensional torus $S$ 
 such that $\Face(C)=\Face$.

There exists a generic closed isotropy $H$ of $\Face$ such that $H^\circ=S$.
\end{prop}

\begin{proof}
By Lemma~\ref{lem:FaceinterC}, $\Face$ is an union of GIT-classes.
By \cite{GeomDedic}, there are only finitely many such classes and they are convex;
so, they exists a GIT-class $F$ which spans $\Face$.
Let  $\Li\in F$.
Let $\xi\in (X^{\rm ss}(\Li)\quot G)_{\rm pr}$ and $\Orb_\xi$ be the corresponding closed $G$-orbit
in $X^{\rm ss}(\Li)$.

By \cite[Theorem~5]{GITEigen}, $\Orb_\xi$ intersects $C$.
Let us fix $x\in O_\xi\cap C$.
Now, consider the morphism $\mu^\bullet(x,G_x)\,:\,\Lambda_\QQ\longto X(G_x)\otimes \QQ$ induced by restriction and 
the isomorphism $X(G_x)\simeq {\rm Pic}^G(\Orb_\xi)$.
By Theorem~\ref{th:CGXloc}, $\Ker\mu^\bullet(x,G_x)$ is contained in $\Span(\Face)$.
On the other hand, the GIT-class of $\Li$ is contained in $\Ker\mu^\bullet(x,G_x)$.
Finally,  $\Ker\mu^\bullet(x,G_x)=\Span(\Face)$.
Since $\Lambda$ is abundant, this implies that the rank of $X(G_x)$ equals $r$.

Since $G.x$ is affine, $G_x$ is reductive. Since $X=Y\times G/B$, $G_x$ is contained in a Borel subgroup of $G$.
Finally, $G_x$ is diagonalisable. But the rank of $X(G_x)$ equals $r$; and, $G_x^\circ$ is a $r$-dimensional torus.

Since $x\in C$, $S$ is contained in $G_x$; it follows that $G_x^\circ=S$.
\end{proof}

\subsection{Unicity}

\begin{nota}
  Let $S$ be a torus. We will denote $Y(S)$ the group of one parameter subgroups of $S$.
There is a natural perfect paring $Y(S)\times X(S)\longto \ZZ$ denoted by $\langle\cdot,\cdot\rangle$. 
\end{nota}

The following proposition is a first statement of unicity. 

\begin{lemma}
\label{lem:unicityCS}
  We assume that   $Y$ (and so $X$) is smooth.
Let $\Face$ be a face of codimension $r$.
Let $(C_1,S_1)$ and $(C_2,S_2)$ be  two  well covering pairs with  two $r$-dimensional tori $S_1$ and 
$S_2$ such that $\Face(C_1)=\Face(C_2)=\Face$.

Then, the there exists $g\in G$ such that $g.C_2=C_1$ and $g.S_2.g^{-1}=S_1$.
\end{lemma}

\begin{proof}
Starting the proof as Proposition~\ref{prop:gciGB}, we obtain that $\Orb_\xi$ intersects $C_1$ and $C_2$.
Up to conjugacy, we may assume that $x$ belongs to $\Orb_\xi\cap C_1\cap C_2$.
So, we obtain that $G_x^\circ=S_1=S_2$.
Then, $C_1$ equals $C_2$ since they are the irreducible component of $X^{S_1}=X^{S_2}$ containing $x$.
\end{proof}\\

Let us fix a face  $\Face$ of codimension $r$. 
The set of linear forms $\varphi\in\Hom(\Lambda_\QQ,\QQ)$ such that $\varphi$ is non negative on $\CGLX$ and
zero on $\Face$ is denote by $\Face^\vee$.

Let $(C,S)$ be an admissible pair where $S$ has dimension $r$ and set
$\Face=\Face(C)$.
Let $\cone$ denote the set of $\lambda\in Y(S)\otimes \QQ$ such that for some positive $n$, the pair 
$(C_{n\lambda},n\lambda)$ is  dominant; where, $C_{n\lambda}$ denote the irreducible component of $X^{n\lambda}$ 
containing $C$.

\begin{lemma}
\label{lem:Facedual}
We assume that $Y$ is smooth.
 
Then, $\lambda\in\cone$ if and only if $\mu^\bullet(C,\lambda)\in\Face^\vee$.   
\end{lemma}

\begin{proof}
Let $\lambda\in\cone$ and $n$ be a positive integer such that $n\lambda\in Y(S)$.
Since $(C_{n\lambda},n\lambda)$ is dominant, \cite[Lemma~7]{GITEigen} imply that $\mu^\bullet(C,\lambda)$ is non negative 
on $\CGLX$. Moreover, for any $\Li\in\Face$, $X^{\rm ss}(\Li)$ intersects $C$. This implies that
$\mu^\Li(C,\lambda)=0$. Finally, $\mu^\bullet(C,\lambda)\in\Face^\vee$.

Conversely, let $\lambda$ be a rational one parameter subgroup and $n$ be a positive integer such that
$n\lambda\in Y(S)$ and $\mu^\bullet(C,\lambda)\in\Face^\vee$.   
Set $C^+=\{x\in X\,:\,\lim_{t\to 0} n\lambda(t)x\in C_{n\lambda}\}$ and
$\eta\,:\,G\times_{P(n\lambda)}C^+\longto X$. 
Let us fix a generic point $x\in C$.
Then, $G_x$ is the generic closed isotropy of $\Face$, its neutral component is $S$ and the $G_x$-module 
$T_xX/T_xG.x$ is the type of $\Face$.
Theorem~\ref{th:CGXloc} implies that $\mu^\bullet(C,\lambda)\in\Face^\vee$ if and only if 
$\langle n\lambda,\cdot\rangle$ is non negative on all weights of $S$ in $T_xX/T_xG.x$.
We deduce that $T\eta_x$ is surjective. Since $Y$ is smooth, this implies that $\eta$ is dominant.
\end{proof}\\

We can now state our main result of unicity:

\begin{prop}
\label{prop:unicity}
We assume that   $Y$ is smooth.
Let $p_{G/B}\,:\,X\longto G/B$ denote the projection. Let us fix a Borel subgroup $B$ of $G$ and a maximal
torus $T$ of $B$.

Let $\Face$ be a face of codimension $r$.
Then there exists a unique well covering pair $(C,S)$ where $S$ is a $r$-dimensional subtorus of $T$,
$p_{G/B}(C)$ contains $B/B$ and $\Face(C)=\Face$.

Let $\lambda_1$ and $\lambda_2$ be two one parameter subgroups of $S$ such that $(C,\lambda_i)$ is dominant and
$\Face$ equals the set of $\Li\in\CGLX$ such that $\mu^\Li(C,\lambda_i)=0$ for $i=1,2$.
Then, $P(\lambda_1)=P(\lambda_2)$, $C$ is an irreducible component of $X^{\lambda_1}$ and $X^{\lambda_2}$ and
$C^+(\lambda_1)=C^+(\lambda_2)$.
\end{prop}

\begin{proof}
Let $(C_1,S_1)$ and $(C_2,S_2)$ be  two admissible well covering pairs with  two $r$-dimensional 
tori $S_1$ and $S_2$ such that $\Face(C_1)=\Face(C_2)=\Face$.
We also assume that $p_{G/B}(C_1)$ and $p_{G/B}(C_2)$ contain $B/B$.
By Lemma~\ref{lem:unicityCS}, there exists $g\in G$ such that $gS_2g^{-1}=S_1$ and $gC_2=C_1$.
Note that $g^{-1}Tg$ and $T$ contain $S_2$ and are maximal tori of $G^{S_2}$:
there exists $h\in G^{S_2}$ such that $hTh^{-1}=g^{-1}Tg$.
Set $\tilde{w}=gh$. One easily checks that $\tilde{w}$ normalize $T$, $\tilde{w}S_2\tilde{w}^{-1}=S_1$ and
$\tilde{w}C_2=C_1$.
Now, $G^{S_1}B/B=p_{G/B}(C_1)=\tilde{w}p_{G/B}(C_2)\ni \tilde{w}B/B$. We deduce that 
$\tilde{w}\in G^{S_1}$.
So, $S_2=\tilde{w}^{-1}S_1\tilde{w}=S_1$ and $C_2=\tilde{w}^{-1}C_1=C_1$.
The first assertion is proved.\\

Let now $C,S,\lambda_1$ and $\lambda_2$ be as in the statement.
Let us first prove that $P(\lambda_1)=P(\lambda_2)$.
By Lemma~\ref{lem:Facedual}, the set of the one parameter subgroups $\lambda$ of $S$ as in the proposition is convex.
So, if it is possible to have $P(\lambda_1)\neq P(\lambda_2)$; it is possible to have $P(\lambda_1)\subset P(\lambda_2)$ 
and $P(\lambda_1)\neq P(\lambda_2)$. In other words, we may assume that $P(\lambda_1)\subset P(\lambda_2)$.
Let $C_1$ (resp. $C_2$) denote the irreducible component of $X^{\lambda_1}$ (resp. $X^{\lambda_2}$) containing $C$.
By Lemma~\ref{lem:unicityCS}, we have $C_1=C_2$. 
In particular, $G^{\lambda_1}B/B=p_{G/B}(C_1)=p_{G/B}(C_2)=G^{\lambda_2}B/B$; so, $G^{\lambda_1}=G^{\lambda_2}$.
It follows that $P(\lambda_1)=P(\lambda_2)$. Let $P$ denote this parabolic subgroup of $G$.

Let $x\in C$ be general.
Since $\lambda_1$ fixes $x$, it acts on the tangent space $T_xX$.
Consider the subspaces $(T_xX)_{>0}$ and $(T_xX)_{<0}$ of the $\xi\in T_xX$ such that 
$\lim_{t\to 0}\lambda_1(t)\xi=0$
and  $\lim_{t\to 0}\lambda_1(t^{-1})\xi=0$ respectively.
We have: $T_xX=(T_xX)_{<0}\oplus (T_xX)_0\oplus (T_xX)_{>0}$. 
The first part identify with $T_eG/P$ and the second one with $T_xC$ as $S$-modules.
But, the third part is the unique $S$-stable supplementary of the sum of the two first one.
In particular, the same construction with $\lambda_2$ in place of $\lambda_1$ gives the same decomposition
 $T_xX=(T_xX)_{<0}\oplus (T_xX)_0\oplus (T_xX)_{>0}$. 
It follows that $C^+(\lambda_1)=C^+(\lambda_2)$.
\end{proof}

\subsection{Inclusion of faces}

\begin{prop}
\label{prop:inclusionface}
  We assume that   $Y$ is smooth.
Let us fix a maximal torus $T$ of $G$ and a Borel subgroup $B$ containing $T$.
Let $(C_1,S_1)$ and $(C_2,S_2)$ be  two admissible well covering pairs with two subtori $S_1$ and 
$S_2$ of $T$ of dimension $r_1$ and $r_2$ such that $B/B\in p_{G/B}(C_i)$ for $i=1,2$.
We assume that $\Face(C_1)$ and $\Face(C_2)$ have respectively codimension $r_1$ and $r_2$.

Then, the following are equivalent:
\begin{enumerate}
\item $\Face(C_1)\subset \Face(C_2)$;
\item $C_1\subset C_2$ and $S_2\subset S_1$.
\end{enumerate}
\end{prop}

\begin{proof}
  The second assertion implies the first one by Lemma~\ref{lem:FaceinterC}.
Conversely, let us assume that  $\Face(C_1)\subset \Face(C_2)$.

By Proposition~\ref{prop:gciGB}, there exists $\Li\in\Face(C_1)$ and $x\in C_1$ such that $G_x^\circ=S_1$ and
$G.x$ is closed in $X^{\rm ss}(\Li)$.
Since $C_2$ intersects $G.x$, there exists $g\in G$ such that $g.x\in C_2$.
Since $S_2$ fixes $g.x$, $S_2\subset gS_1g^{-1}$.
In particular, $S_2$ is contained in $T$ and $gTg^{-1}$; so, 
$T$ and $gTg^{-1}$ are maximal tori in $G^{S_2}$.
There exists $l\in G^{S_2}$ such that $lTl^{-1}=gTg^{-1}$.
The element  $n=l^{-1}g$ normalizes $T$.
Since $C_2$ is stable by $G^{S_2}$ (which is connected), 
$x$ belongs to $n^{-1}C_2$.
Applying $p_{G/B}$ we obtain that $p_{G/B}(x)$ belongs to $n^{-1}G^{S_2}B/B\cap G^{S_1}B/B$.

Since $n^{-1}S_2n=g^{-1}S_2g\subset S_1$, we have $G^{n^{-1}S_2n}\subset G^{S_1}$. 
In particular,  $n^{-1}G^{S_2}B/B\subset G^{S_1}n^{-1}B/B$.
It follows that $G^{S_1}n^{-1}B/B=G^{S_1}B/B$.

Since $n$ normalizes $T$, this implies that $n\in G^{S_1}$.
So, $S_2\subset nS_1n^{-1}=S_1$.\\

Since $n\in G^{S_1}$, $nx\in C_1$.
It follows that $C_1$ is  the irreducible component of $X^{S_1}$ containing $nx$.
On the other hand, $C_2$ is  the irreducible component of $X^{S_2}$ containing $nx$.
It follows that $C_1$ is contained in $C_2$.
\end{proof}

\section{GIT-cone and moment polytope}
\label{sec:momentpol}

Let us now explain the relation mentioned in the introduction between 
the moment polytopes of $Y$ and some total $G$-cones of $X=Y\times G/B$.

Let us fix a maximal torus $T$ of $G$ and a Borel subgroup $B$ containing $T$.
Let $\Li$ be an ample $G$-linearized line bundle on $Y$. We consider the set $P_G(Y,\Li)$ of the points 
$p\in\Chi(T)_\QQ$ such that for some positive integer $n$, $np$ is a dominant character of $T$ and 
the dual $V_{np}^*$ of $V_{np}$ 
is a submodule of ${\rm H}^0(Y,\Li^{\otimes n})$. 
In fact, $P_G(Y,\Li)$ is a polytope, called {\it moment polytope}.
Notice that ``the dual'' is not usual in the definition; but it will be practical in our context. 

Consider the two projections:
$$
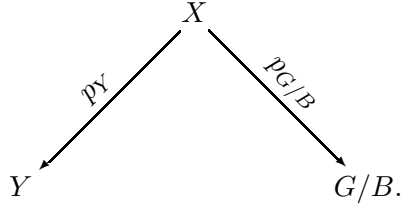
\begin{diagram}
  &&X&&\\
&\ldTo^{p_Y}&&\rdTo^{p_{G/B}}&\\
Y&&&&G/B.
\end{diagram}
$$
In Section~\ref{sec:momentpol}, $\Lambda$ will always denote the subgroup of $\Pic^G(X)$ generated by 
$p_{G/B}^*(\Pic^G(G/B))$ and  $p_Y^*(\Li)$.
Note that $\Pic^G(G/B)$ identifies with $X(T)$; we will denote $\Li_\nu$ the element of $\Pic^G(G/B)$ 
associated to $\nu\in X(T)$

\begin{prop}
With above notation, we have:
\begin{enumerate}
\item $\tc^G_\Lambda(X)=\sac^G_\Lambda(X)$;
\item $\sac^G_\Lambda(X)$ is a cone over $P_G(Y,\Li)$;
more precisely, for all positive rational number $m$ and $\nu\in\Chi(T)$, we have:
$$
mp_Y^*(\Li)\otimes p_{G/B}^*(\Li_\nu)\in\sac^G_\Lambda(X)\iff
\frac{\nu}{m}\in P_G(Y,\Li).
$$
\end{enumerate}
\end{prop}

\begin{proof}
Let $n$ be a non negative integer and $\nu$ be a character of $T$.
As a $G$-module, ${\rm H}^0(X,p_Y^*(\Li^{\otimes n})\otimes p_{G/B}^*(\Li_\nu))$ is isomorphic
to ${\rm H}^0(Y,\Li^{\otimes n})\otimes {\rm H}^0(G/B,\Li_\nu)$.
In particular, if $p_Y^*(\Li^{\otimes n})\otimes p_{G/B}^*(\Li_\nu)$ has non zero global sections then
$n\geq 0$ and $\nu$ is dominant; in this case, it is semiample. The first assertion follows.

Assume now that $\nu$ is dominant. Then, ${\rm H}^0(G/B,\Li_\nu)=V_\nu$.
Hence,  $p_Y^*(\Li^{\otimes n})\otimes p_{G/B}^*(\Li_\nu)$ has non zero $G$-invariant sections if and only if 
 ${\rm H}^0(Y,\Li^{\otimes n})\otimes V_\nu$ contains non zero $G$-invariant vectors; that is, if 
and only if $V_\nu^*$ is a submodule of ${\rm H}^0(Y,\Li^{\otimes n})$. The second assertion follows.
\end{proof}\\

\begin{remark}
To each face $\Face$ of $\CGLX$, one can associate a face of $P_G(Y,\Li)$ (by intersecting and taking a closure) 
which intersects the interior of the dominant chamber. 
By this way, we obtain a bijection between the set of faces of $\CGLX$ and the faces of $P_G(Y,\Li)$ which   
intersects the interior of the dominant chamber. 
\end{remark}

\subsection{A reduction}

It is possible that $\ac_\Lambda^G(X)$ is empty. In this case, our results cannot be
applied directly.

Let $mp_Y^*(\Li)\otimes p_{G/B}^*(\Li_\nu)$ in the relative interior of $\sac^G(X)$ such that 
${\rm  H}^0(X,mp_Y^*(\Li)\otimes p_{G/B}^*(\Li_\nu))^G$ is non zero, that is
such that $V_\nu^*$ can be equivariantly embedded in ${\rm  H}^0(Y,\Li^{\otimes m})$.
Let $P$ be the standard parabolic subgroup of $G$ associated to the face of the dominant chamber containing $\nu$.
Let $L$ denote the Levi subgroup of $P$ containing $T$ and $D$ denote its derived subgroup.

The next proposition shows that $\sac^G(X)$ is always equal to such a cone satisfying $\ac^G(X)^neq\emptyset$.
The proof which is essentially extracted from \cite[Section~5]{Br:genface} is included for completeness.

\begin{prop}
\label{prop:PXLempty}
With above notation, there exists an irreducible component $C_Y$ of $Y^D$ such that 
a point $\Li\in\Lambda_\QQ$ belongs to $\sac^G(X)$ if and only if 
$\Li_{|C}$ belongs to $\sac^{G^D}(C_Y\times G^D.B/B)$.

Moreover, $G^D.B/B$ is isomorphic to the variety of complete flags of $G^D$ and
$\ac^{G^D}(C_Y\times G^D.B/B)$ is non empty.
\end{prop}

\begin{proof}
  The inclusion $V_\nu^*\subset{\rm  H}^0(Y,\Li^{\otimes m})$ gives a $G$-equivariant rational map
$\phi\::\:Y\longratto\PP(V_\nu)$.
Let $v_\nu$ be a vector of highest weight in $V_\nu$;
$P$ is the stabilizer in $G$ of $[v_\nu]\in\PP(V_\nu)$. 
Let $\sigma\in V_\nu^*$ be an eigenvector of the Borel subgroup $B^-$ opposite to $B$ and containing $T$.  
Let $Q$ be the stabilizer in $G$ of $[\sigma]\in\PP(V_\nu^*)$ and $Q^u$ be its unipotent radical.

Let $Y_\sigma$ denote the set of $y\in Y$ such that $\sigma(y)\neq 0$.
Let $W$ be a $L$-stable supplementary subspace of $\kk.v_\nu$ in $V_\nu$.
By $w\longmapsto [v_\nu+w]$, we identify $W$ with an open subspace of $\PP(V_\nu)$.
Then $\phi$ induces by restriction $\widetilde{\phi}\::\:Y_\sigma\longto W$. 

Let $S$ be a $L$-stable supplementary to $T_{[v_\nu]}G.[v_\nu]$ in $W$ 
(actually, $W$ canonically identify with $T_{[v_\nu]}\PP(V_\nu)$).
Set $Z=\widetilde{\phi}^{-1}(S)$.
By \cite[Remark in Section~5]{Br:genface}, $Z$ is point wise fixed by $D$ and 
the action of $Q^u$ induces an isomorphism $Q^u\times Z\simeq Y_\sigma$.\\

Consider $X'=Y\times G/P$.
Let $\Lambda'$ be the subgroup of $\Pic^G(X')$ generated by $p_Y^*(\Li)$ and
$p_{G/P}^*(\Pic^G(G/P))$ (with obvious notation). 
It is clear that $\tc^G_\Lambda(X)$ identifies with  $\tc^G_{\Lambda'}(X')$.
Moreover, $\ac^G_{\Lambda'}(X')$ is not empty.
Consider a generic closed isotropy $H$ of $\tc^G_{\Lambda'}(X')$ viewed as a face $\Face$ of itself.
Since  $Q^u\times Z\simeq Y_\sigma\subset X^{\rm ss}(mp_Y^*(\Li)\otimes p_{G/P}^*(\Li_\nu))$, up to conjugacy, 
one may assume that $D\subset H\subset L$.
Since $Y_\sigma\times \{P/P\}\subset  X^{\rm ss}(mp_Y^*(\Li)\otimes p_{G/P}^*(\Li_\nu))$,
$X'_\Face$ intersects $Z\times \{P/P\}$. 

Consider the irreducible component $C_Y$ of $Y^D$ which contains $Z$.
By Theorem~\ref{th:HFCF}, for any ample $\Li\in\Lambda'$, 
$X'^{\rm ss}(\Li)$  intersects $C_Y\times \{P/P\}$ if it is non empty.
By continuity, this is also true if $\Li$ is only semiample. 
The proposition follows easily.
\end{proof}

\subsection{Faces of $\sac^G(X)$ if $\ac^G(Y)$ is non empty}

From now on, we assume that $Y$ is smooth.
We will first adapt the notion of covering and well covering pairs for the
situation.

Recall that $T\subset B$ are fixed.
Let $\lambda$ be a one parameter subgroup of $T$. Set $B(\lambda)=B\cap P(\lambda)$.
Let $C$ be an irreducible component of $Y^\lambda$ and $C^+$ the associated Bialinicki-Birula cell.\\

\begin{defin}
\label{def:wellBcov}
  The pair $(C,\lambda)$ is said to be {\it $B$-covering} if the natural map
$\eta\::\:B\times_{B(\lambda)}C^+\longto Y$ is birational. 
It is said to be {\it well $B$-covering} is $\eta$ induces an isomorphism over an open subset of $Y$ 
intersecting $C$.
\end{defin}

The proof of the following lemma is obvious.

\begin{lemma}
\label{lem:covBcov}
With above notation, the pair $(C,\lambda)$ is $B$-covering (resp. well $B$-covering) if and only if
$(C\times G^\lambda B/B)$ is covering (resp. well covering).   
\end{lemma}

Let us recall that the subtori of $T$ correspond bijectively to the linear subspaces of $\Chi(T)_\QQ$.
If $V$ is a linear subspace of $\Chi(T)_\QQ$, the associated torus is the neutral component of the 
intersection of  kernels of elements in $\Chi(T)\cap V$.
If $F$ is a convex part of $\Chi(T)_\QQ$, the direction ${\rm dir}(F)$ of $F$ is the linear subspace
 spanned by the differences of two elements of $F$.

We will denote by $\cone^+$ the convex cone in $X(T)_\QQ$ generated by the dominant weights.
The next proposition is an improvement of \cite[Theorem~1]{Br:genface}:

\begin{prop}
\label{prop:PYL}
We keep the above notation and assume that $Y$ is smooth and $P_G(Y,\Li)$ intersects the
interior of the dominant chamber.
Let $\Face$ be a face of codimension $d$ of $P_G(Y,\Li)$ which intersect the interior of the 
dominant chamber.
Let $S$ the subtorus of $T$ associated to ${\rm dir}(\Face)$.

There exists a unique irreducible component $C$ of $Y^S$ and a one parameter subgroup $\lambda$ 
of $S$ such that $G^\lambda=G^S$ and $(C,\lambda)$ is a well $B$-covering pair such that 
$\Face=P_{G^S}(C,\Li_{|C})\cap \cone^+$.
\end{prop}

\begin{proof}
Let $\widetilde{\Face}$ be the face of $\CGLX$ corresponding to $\Face$ and $r$ denote its codimension.
By Proposition~\ref{prop:ttFC}, there exists an admissible well covering pair $(C_X,S')$ such that 
$\widetilde{\Face}=\Face(C_X)$ and $S'$ is a $r$-dimensional torus.
Up to conjugacy, we may assume that $C_X$ intersects $Y\times B/B$, and $S'$ is contained in $T$.
Let $\lambda$ be a one parameter subgroup of $S'$ such that $(C_X,\lambda)$ is well covering.
Then, $C_X=C\times G^\lambda B/B$ for some irreducible component $C$ of $Y^{S'}$.

The fact  $\widetilde{\Face}=\Face(C_X)$ readily means that $\Face=P_{G^{S'}}(C,\Li_{|C})\cap \cone^+$.
Since the direction of $P_{G^{S'}}(C,\Li_{|C})$ is contained in $X(T)^{S'}$, this implies that $X(T)^S$ is
contained in $X(T)^{S'}$.
But, $S$ and $S'$ have the same rank, it follows that $S=S'$.

The unicity part is a direct consequence of Proposition~\ref{prop:unicity}.
\end{proof}

\section{The case $X=\hG/\hB\times G/B$}
\label{sec:BeSj}

\subsection{Interpretations of the $G$-cones}

From now on, we assume that $G$ is a connected reductive subgroup of a connected reductive group $\hG$.
Let us fix  maximal tori $T$ (resp. $\hT$) and Borel subgroups $B$ (resp. $\hB$) of $G$ (resp. $\hG$) such that 
$T\subset B\subset \hB\supset\hT\supset T$.

Let $\lg$ and $\hlg$ denote the Lie algebras of $G$ and $\hG$ respectively.

We denote by $\lr(\hG,G)$ (resp. $\lr^\circ(\hG,G)$) the cone of the pairs 
$(\hnu,\nu)\in \Chi(\hT)_\QQ\times\Chi(T)_\QQ$ 
such that for a positive integer $n$, $n\hnu$ and $n\nu$ are dominant (resp. strictly dominant) 
weights such that $V_{n\hnu}\otimes V_{n\nu}$ contains non zero $G$-invariant vectors.

In this section, $X$ denote the variety $\hG/\hB\times G/B$ endowed with the diagonal action of $G$.
We will apply the results of Section~\ref{sec:gen} to $X$ with $\Lambda=\Pic^G(X)$.
The cones $\tc^G(X),\,\sac^G(X)$ and $\ac^G(X)$ will be denoted without the $\Lambda$ in subscribe.
By \cite[Proposition~9]{GITEigen}, $\lr^\circ(\hG,G)=\ac^G(X)\subset\sac^G(X)=\tc^G(X)=\lr(\hG,G)$.
Moreover, if no ideal of $\lg$ is an ideal of $\hlg$, 
by \cite[Assertion~$(i)$ of Theorem~9]{GITEigen} $\lr^\circ(\hG,G)$ has 
non empty interior.

\subsection{The case $X=\hG/\hB\times G/B$}
\label{sec:GBGB}

\paragraphe
Consider the $G$-module $\hlg/\lg$. 
Let $\chi_1,\cdots,\chi_n$ be the set of the non trivial weights of $T$ on $\hlg/\lg$.
For $I\subset\{1,\cdots,n\}$, we will denote by $T_I$ the neutral component of the intersection
of the kernels of the $\chi_i$'s with $i\in I$. 
A subtorus of the form $T_I$ is said to be {\it admissible}.

Let $\lambda$ be a one parameter subgroup of $T$.
Consider the parabolic subgroups $P$ and $\hP$ of $G$ and $\hG$ associated to $\lambda$.
Let $W_P$ be the Weyl group of $P$.
The cohomology group ${\rm H}^*(G/P,\ZZ)$ is freely generated by the Schubert classes
$[\overline{BwP/P}]$ parametrized by the  elements $w\in W/W_P$.
Since $\hP\cap G=P$, we have a canonical $G$-equivariant immersion 
$\iota\,:\,G/P(\lambda)\longto\hG/\hP(\lambda)$; and the corresponding morphism 
in cohomology $\iota^*$. 

Let $\rho$ (resp. $\rho^\lambda$) 
denote the half sum of the positive roots of $G$ (resp. $G^\lambda$).
Let $\Phi^+$ and $\Phi(P^u)$ denote the set of roots of the groups
$B$ and $P^u$ for the torus $T$. In the same way, we define $\hat\Phi^+$ and $\Phi(\hP^u)$.
For $\hw\in \hW$, we set:
$$
\theta^P:=\sum_{\alpha\in\Phi^+\cap \Phi(P^u)}\alpha\in X(T)\ \ \ {\rm and}\ \ \ 
 \theta^\hP_\hw:=\sum_{\alpha\in\hw\hat\Phi^+\cap \Phi(\hP^u)}\alpha\in X(\hT).
$$

\paragraphe
Let $S$ be an admissible subtorus of $T$.
All irreducible component $C$ of $X^S$ such that $p_{G/B}(C)$ contains $B/B$ equals
$C(\hw):=(\hG^S.\hw^{-1}\hB/\hB\times G^SB/B)$ for a unique element $\hw\in\hW/\hW_{\hG^S}$. 
Let us fix $\hat w\in \hW/\hW_{\hG^S}$.
The pair $(S,\hw)$ is said to be {\it admissible} if there exists a parabolic subgroup $\hP$ of $\hG$
such that
\begin{enumerate}
\item there exists $\lambda\in Y(S)$ such that $\hP=\hP(\lambda)$;
\item $\hG^S$ is a Levi subgroup of $\hP$;
\item $G^S$ is a Levi subgroup of $\hP\cap G=:P$;
\item $\iota^*([\overline{\hB\hw\hP/\hP}]).[\overline{BP/P}]=[{\rm pt}]\in{\rm H}^*(G/P,\ZZ)$;
\item $(\theta^\hP_\hw)_{|S}=(\theta^P-2(\rho-\rho^\lambda))_{|S}$.
\end{enumerate}

\begin{lemma}
\label{lem:wadmissible}
Let $S$ be an admissible subtorus of $T$ and $\hat w\in \hW/\hW_{\hG^S}$.
The pair $(S,\hw)$ is admissible if and only if there exists a one parameter subgroup $\lambda$ of $S$
such that $C(\hw)$ is an irreducible component of $X^\lambda$ and $(C(\hw),\lambda)$ is a well covering pair.
\end{lemma}

\begin{proof}
  The proof is very analogous to \cite[Proposition~10]{GITEigen}: we leave details to the reader.
We prove (using mainly Kleiman's Theorem) that 
$\iota^*([\overline{\hB\hw\hP/\hP}]).[\overline{BP/P}]=[{\rm pt}]\in{\rm H}^*(G/P,\ZZ)$
if and only if $\eta$ is birational.
Now, the condition  $(\theta^\hP_\hw)_{|S}=(\theta^P-2(\rho-\rho^\lambda))_{|S}$ means that $S$ 
acts trivially on the restriction over $C$ of the determinant bundle of $\eta$. 
\end{proof}\\

\paragraphe
To simplify, in the following statement we assume that $\ac^G(X)$ has a non empty interior
in $\Pic^G(X)_\QQ$. In fact, this assumption is equivalent to say that no ideal of $\lg$ is an 
ideal of $\hlg$.

\begin{theo}
\label{th:ppal}
  We assume that no ideal of $\lg$ is an ideal of $\hlg$.

The map which associates to a pair $(S,\hw)$ the set 
$\Face(S,\hw)=\{(\hnu,\nu)\in C^G(X)\,:\,\hw\hnu_{|S}=-\nu_{|S}\}$ is a bijection from the set 
of admissible pairs onto the set of faces of $\CGLX$.
Moreover, the codimension of $\Face(S,\hw)$ equals the dimension of $S$.

The following are equivalent:
\begin{enumerate}
\item $\Face(S,\hw)\subset\Face(S',\hw')$;
\item $S'\subset S$ and $\hw W_{G^{S'}}=\hw' W_{G^{S'}}$.
\end{enumerate}
\end{theo}

\begin{proof}
Let $(S,\hw)$ be an admissible pair. 
Set $\overline{\Face}(S,\hw)=\{(\hnu,\nu)\in \lr(G,\hG)\,:\,\hw\hnu_{|S}=-\nu_{|S}\}$.
By \cite[Theorem~9]{GITEigen},  $\overline{\Face}(S,\hw)$ is  a face of $\lr(G,\hG)$ 
of codimension ${\rm dim}(S)$.
In particular, $\widetilde{\Face}(S,\hw)$ spans the vector subspace of the 
$(\hnu,\nu)\in X(\hT)\times X(T)$ such that $\hw\hnu_{|S}=-\nu_{|S}\}$.
Since $\CGLX$ is the interior of $\lr(G,\hG)$, to prove that the map in the theorem is well defined, 
it is enough to prove that $\overline{\Face}(S,\hw)$ intersects $\CGLX$.
If not,  $\overline{\Face}(S,\hw)$ would be contained in the boundary of the dominant chamber.
Its projection on $X(\hT)_\QQ$ or $X(T)_\QQ$ would be contained in an hyperplane;
which is a contradiction. 

The surjectivity is a rephrasing of \cite[Theorem~9 Assertion $(ii)$]{GITEigen}.
The injectivity is a direct application of Proposition~\ref{prop:unicity}.

The last assertion follows from Proposition~\ref{prop:inclusionface}.
\end{proof}\\

\section{Application to the tensor product cone}

\paragraphe
In this section, $G$ is assumed to be semisimple.
As above, $T\subset B$ are fixed maximal torus and Borel subgroup of $G$. 
We also fix an integer $s\geq 2$ and set $\hG=G^s$, $\hT=T^s$ and $\hB=B^s$.
We embed $G$ diagonally in $\hG$. 
Now, $X=\hG/_hB\times G/B=(G/B)^{s+1}$.
Then $\ac^G(X)\cap \Chi(T)^{s+1}$ identifies with the $(s+1)$-uple 
$(\nu_1,\cdots,\nu_{s+1})\in\Chi(T)^{s+1}$ such that the for $n$ big enough $n\nu_i$'s are 
strictly dominant weights and $V_{n\nu_1}\otimes\cdots\otimes V_{n\nu_{s+1}}$ contains a non zero 
$G$-invariant vector.

A parabolic subgroup $P$ of $G$ is said to be {\it standard} if it contains $B$. 
We will denote by $Z(P)$ the neutral component of the center of the Levi subgroup of $P$ 
containing $T$.

\paragraphe
In \cite{BK}, Belkale and Kumar defined a new product denoted $\kbprod$ on the cohomology
groups $H^*(G/P,\ZZ)$ for any parabolic subgroup $P$ of $G$.
We consider the set $\Theta$ of the $(P,\Lambda_{w_0},\cdots,\Lambda_{w_s})$ where $P$ is a standard 
parabolic subgroup of $G$ and the $\Lambda_{w_i}$'s are $s+1$ Schubert varieties of $G/P$ such that
$$
[\Lambda_{w_0}]\kbprod\cdots\kbprod[\Lambda_{w_s}]=[{\rm pt}].
$$

\paragraphe
In \cite{GITEigen}, Theorem~9 applied to $\hG=G^s$ gives Corollary~5.
The same translation of Theorem~\ref{th:ppal} to this case gives the following:

\begin{theo}
The map which associates to a  $(P,\Lambda_{w_0},\cdots,\Lambda_{w_s})\in \Theta$ the set 
$\Face(P,\Lambda_{w_0},\cdots,\Lambda_{w_s})$ of the 
$(\nu_0,\cdots,\nu_s)\in \ac^G(X)$ such that the restriction of $\sum_iw_i^{-1}\nu_i$
to $Z(P)$ is trivial is a bijection from $\Theta$ onto the set of faces of $\ac^G(X)$.
Moreover, the codimension of $\Face(P,\Lambda_{w_0},\cdots,\Lambda_{w_s})$ equals the dimension 
of $Z(P)$.

The following are equivalent:
\begin{enumerate}
\item $\Face(P,\Lambda_{w_0},\cdots,\Lambda_{w_s})\subset\Face(P',\Lambda_{w'_0},\cdots,\Lambda_{w'_s}))$;
\item $P\subset P'$ and $\pi(\Lambda_{w_i})=\Lambda_{w'_i}$ for all $i=0,\cdots,s$ (here,
$\pi\,:\,G/P\longto G/P'$ is the natural $G$-equivariant map).
\end{enumerate}  
\end{theo}

\bibliographystyle{amsalpha}
\bibliography{biblio}

\begin{center}
  -\hspace{1em}$\diamondsuit$\hspace{1em}-
\end{center}

\vspace{5mm}
\begin{flushleft}
N. R.\\
Universit{\'e} Montpellier II\\
D{\'e}partement de Math{\'e}matiques\\
Case courrier 051-Place Eug{\`e}ne Bataillon\\
34095 Montpellier Cedex 5\\
France\\
e-mail:~{\tt ressayre@math.univ-montp2.fr}  
\end{flushleft}

\end{document}